\documentclass[12pt]{article}
\usepackage{amsmath,amsthm,amssymb}
\usepackage{german}
\theoremstyle{plain} %% This is the default
\newtheorem{thm}{Theorem}
\newtheorem{cor}[thm]{Corollary}
\newtheorem{lem}[thm]{Lemma}

\theoremstyle{definition}
\newtheorem{defn}{Definition}

\newcommand{\ncr}[2]{\mbox{$\left(\begin{array}{c}#1\\#2\end{array}\right)$}}

\title{
A reconstruction problem related to balance equations-I
\footnote{Published in Discrete Mathematics  176 (1997) 279-284.}
}
\author{Bhalchandra D. Thatte \\
\small \texttt{bdthatte@gmail.com} 
}
\date{}

\begin{document}
\maketitle

%%fakesection abstract
\begin{abstract}
A modified $k$-deck of a graph is obtained by removing $k$ edges in all 
possible ways and adding $k$ (not necessarily new) edges in all 
possible ways. Krasikov and Roditty used these decks to give an 
independent proof of 
M\"uller's result on the edge reconstructibility of graphs.
They asked if a   
$k$-edge deck could be constructed from its modified
$k$-deck. In this paper, we solve the problem when $k=1$. 
We also offer new proofs  
of Lov\'asz's result, one describing the constructed graph explicitly,
(thus answering a question of Bondy), and another 
based on the eigenvalues of Johnson graph. 
\end{abstract}

\section{Introduction}
\label{intro}
Graphs in this paper  
are 
assumed to be undirected and without multiedges or loops,  
and are assumed to have $n$ vertices and $m$ edges, unless
specified otherwise. The complement of $G$ is denoted by $G^c$. 
We set  $N = \ncr{n}{2}$. 

For a graph $G$, we define three kinds of decks -- $MD_i(G)$, 
$PD_i(G)$ and $ED_i(G)$ as follows. The deck $MD_i(G)$,
called the modified $i$-deck in [KR], is 
obtained by removing $i$ edges of $G$ in all possible ways 
and adding $i$ (not necessarily new) edges in all 
possible ways.  
When we insist that the set of added edges and the 
set of removed edges be disjoint, we 
call it the perturbed $i$-deck, and denote it  
by $PD_i(G)$. 
The deck $ED_i(G)$ is the $i$-edge deck,  
i.e., the 
collection of all the subgraphs obtained by deleting $i$ edges.
We similarly define the above decks for a collection 
$S=\{G_1,G_2,...,G_r\}$ of graphs as the multiunion of the 
corresponding decks for 
all the members of the collection. Thus, e.g., $MD_i(S) 
= \cup_j MD_i(G_j)$, 
where $\cup $ denotes a multiunion. 

In this paper, we denote sets (or multisets)
by their characteristic vectors of length equal to 
the number of nonisomorphic graphs on $n$ vertices and $m$ edges,
(or $m-i$ edges in case of the $i$-edge deck). The characteristic vector of 
a set $P$ of graphs is denoted by $X_P$, and that of the singleton set
$\{G\}$ by simply $X_G$. The operations 
of constructing the above decks are represented by matrices 
whose rows are indexed by all the nonisomorphic graphs on 
$n$ vertices and $m$ edges, ($m-i$ edges when we are considering 
$i$-edge decks),
and whose columns are indexed by nonisomorphic 
graphs on $n$ vertices and $m$ edges. We denote the matrices corresponding to 
$MD_i$, $PD_i$ and $ED_i$ respectively by  
$\Delta_i$, $D_i$ and $d_i$. The notation chosen here
for denoting the matrices is same as the notation used in [KR]
for the corresponding decks. The $kl$-th entry of $\Delta_i$ is 
the number of graphs isomorphic to $G_k$ in $MD_i(G_l)$.  
The $kl$-th entry of $D_i$ is 
the number of graphs isomorphic to $G_k$ in  $PD_i(G_l)$. 
The $kl$-th entry of $d_i$ is 
the number of graphs isomorphic to $G_k$ in  $ED_i(G_l)$. 
Note that matrices $D_0$ and $\Delta_0$ 
are identity operators. We derive various useful identities involving 
linear combinations of these 
matrices or polynomials in the matrices.

Krasikov and Roditty [KR] used the modified decks to set up their 
balance equations, which they used to derive M\"uller's result independently.
It is easy to see that  $MD_i(G)$ is constructed 
from  $ED_i(G)$ by simply adding $i$ edges in all possible ways in 
each graph in  $ED_i(G)$. 
They asked if  $ED_i(G)$ could be constructed 
from  $MD_i(G)$.  
This problem is solved for the case when $i=1$.  
Also, a new proof of Lov\'asz's result based on the eigenvalues of 
Johnson graph is obtained. 

\section{Reconstructing $ED_1$ from $MD_1$.}

\label{recmod}

We first derive some identities which various matrix operators 
satisfy. 

\begin{lem}
\label{DelD}
$\Delta_k = \sum_{i=0}^k \ncr{m-i}{k-i}D_i $
\end{lem}

\noindent {\bf Proof} This equivalent to Lemma 3.1 in [KR]. 
\begin{thm} 
\label{recursion}
\begin{displaymath}
D_1D_i=(m-i+1)(N-m-i+1)D_{i-1}
+ i(N-2i)D_i + (i+1)^2D_{i+1}
\end{displaymath}
\end{thm} 

\noindent {\bf Proof} For a graph $G$, consider a typical 
member $G - X + Y$ 
of  $PD_i(G)$, where $X\cap Y=\phi$. 
For an edge $e \in E(G) - X + Y$ and 
and an edge $f \not \in E(G) - X + Y$, there is a graph   
$H= G-X+Y-e+f$ in  $PD_1(PD_i(G))$.
Depending upon the choice of $e$ and $f$, 
following four cases arise.
\begin{enumerate}
\item $e \in Y $ and $f\in X $: $H \in PD_{i-1}(G)$. 
\item $e \in E(G)-X$ and $f \in X$: $H \in PD_i(G)$. 
\item $e\in Y$ and $f \in E(G^c)-Y$: 
$H \in PD_i(G)$. 
\item $e\in E(G)-X$ and $f \in E(G^c)-Y$: $H \in PD_{i+1}G$.
\end{enumerate}
In case 1,  same $H$ would be obtained if we first replaced $X'  = X -f +f' $
by $Y'= Y - e + e'$,  where $f\neq f' $ and $e\neq e' $,
and then $e'$ by $f' $.
There are $(m-i+1)$ choices for $f' $ and $(N-m-i+1)$ for
$e'$. This explains the first term on the right hand side.
In a similar manner, the second term results from cases 2 and 3, and the third 
term from
case 4.

Theorem~\ref{recursion} is the basis of our procedure of relating
$\Delta_{i}$ to $\Delta_1$.
We describe the procedure by first relating $\Delta_2$ to $\Delta_1$. 
We apply Theorem~\ref{recursion} for $i=1$. That allows us to express $D_2$
in terms of $D_0$, $D_1$, and $D_1D_1$. 
Since, $D_1$ can be expressed in terms of $\Delta_1 $ and $D_0$,  we 
can express $\Delta_2$ in terms of $D_0$, $\Delta_1$ and $\Delta_1^2$.
It is then enough to demonstrate that the expression for 
 $\Delta_2$ obtained this way has 
no $D_0$ term.  From Lemmas~\ref{DelD} and Theorem~\ref{recursion}, we write 
\begin{eqnarray*}
  \Delta_2 & = & \ncr{m}{2} + \ncr{m-1}{1}D_1 + \ncr{m-2}{0}D_2 \\
  & = & \frac{m(m-1)}{2}D_0 + (m-1)(\Delta_1 - mD_0)\\
  & & +
\frac{1}{4}\{(\Delta_1 - mD_0)(\Delta_1 -  mD_0)-m(N-m)D_0 
- (N-2)(\Delta_1- mD_0)\} \\
  & = & \frac{2(m-1)-N}{4}\Delta_1 + \frac{1}{4}\Delta_1^2
\end{eqnarray*}

To relate $\Delta_i$, $i > 1$ to $\Delta_1$, we first 
prove a lemma  which relates $D_i$ to $\Delta_1$ and $D_0$, and 
then use Lemma~\ref{DelD}.

\begin{lem}
\label{DDel}
$D_k$ can be 
written as \\
 $D_k = L_k + (-1)^k\ncr{m}{k}D_0 $   
where $L_k$ is a linear combination of operators constructed from 
$\Delta_1$.  
\end{lem} 

\noindent {\bf Proof} We prove this by induction on $k$. 
For $k=0$ and $k=1$, the result is trivial to prove. Let 
$D_k = L_k + (-1)^k\ncr{m}{k}D_0 $ for $k \leq r$. 
To prove the result for $k=r+1$, 
we operate from left by $D_1$ on both sides of  
$D_r = L_r + (-1)^r\ncr{m}{r}D_0$.  
From Theorem~\ref{recursion} and the induction hypothesis we get, \\
\begin{eqnarray*} 
& & 
D_1 \left\{ L_r + (-1)^r\ncr{m}{r}D_0 \right\} \\
 & = &  
(m-r+1)(N-m-r+1)\left\{ L_{r-1}+(-1)^{r-1}\ncr{m}{r-1}D_0\right\}  \\
 & & +r(N-2r)\left\{ L_r+ (-1)^r\ncr{m}{r}D_0\right\} 
 + (r+1)^2D_{r+1}
\end{eqnarray*} 
Since $D_1(D_0) = D_1 = \Delta_1-mD_0$, we have
\begin{eqnarray*} 
& & (r+1)^2 D_{r+1} \\
 & = & D_1(L_r)+(-1)^r\ncr{m}{r}\Delta_1 - 
(m-r+1)(N-m-r+1)L_{r-1} \\
& & - r(N-2r)L_r - m (-1)^r\ncr{m}{r}D_0 \\
& & -(m-r+1)(N-m-r+1)(-1)^{r-1}\ncr{m}{r-1}D_0 \\
& & -r(N-2r)(-1)^r\ncr{m}{r}D_0         
\end{eqnarray*}  

From the induction hypothesis, we claim that first 
four terms on the right side give a linear combination of operators 
constructed from $\Delta_1$. Also,  
coefficient of $D_0$ on the right side simplifies to 
$(-1)^{r+1}(r+1)^2\ncr{m}{r+1}$.
Thus $D_{r+1} = L_{r+1}  + (-1)^{r+1}\ncr{m}{r+1}D_0 $, which 
completes the proof.

\begin{thm} 
\label{Delk}
For $i\geq 1$, $\Delta_i $ can be written in terms of 
$\Delta_1 $ 
\end{thm}  

\noindent {\bf Proof} 
In the expression for $\Delta_i$ given by Lemma~\ref{DelD}, 
we substitute    
$D_k $ from Lemma~\ref{DDel},   
Then it is enough to prove that the the coefficient of $D_0$ in 
the resulting expression is 
zero. 
\begin{eqnarray*}
\Delta_i & = & \sum_{k=0}^i\ncr{m-k}{i-k}D_k \\  
 & = &  \sum_{k=0}^i\ncr{m-k}{i-k}\{(-1)^k \ncr{m}{k}D_0+ L_k\} \\
 & = & \sum_{k=0}^i\ncr{m-k}{i-k}(-1)^k \ncr{m}{k}D_0 + \sum_{k=0}^i L_k \\ 
 & = & \sum_{k=0}^i L_k  
\end{eqnarray*}

\begin{lem}
\label{lov}
Let $m \geq N/2$, and let $P$ and $Q$ be collections of graphs 
such that $d_1X_P=d_1X_Q$ or  
$\Delta_1X_P=\Delta_1X_Q$,  then  $X_P=X_Q$.   
\end{lem}

\noindent {\bf Proof} 
Since modified decks can be obtained from edge decks, we assume 
$\Delta_1X_P=\Delta_1X_Q$. 
From Theorem~\ref{Delk}, 
$\Delta_kX_P = \Delta_kX_Q$ for $k\geq 1$. 
Now the claim in the 
lemma is an immediate consequence of the fact that  the
version of Lov\'asz's result presented in
Remark 1 in the next section
makes use of the modified
decks only.
We also note that Krasikov and Roditty used only the 
modified decks in their proof of M\"uller's result. 

Now we prove the main result of this section.
\begin{thm} 
\label{Deld}
For collections $P$ and $Q$ of graphs, if $\Delta_1X_P=\Delta_1X_Q$ then 
$d_1X_P=d_1X_Q$. 
\end{thm} 

\noindent {\bf Proof} 
Let $P' = \{F^c; F\in ED_1(P)\}$ and 
$Q' = \{F^c; F\in ED_1(Q) \}$. 
Note that  $\Delta_1X_P=\Delta_1X_Q$
is equivalent to  $d_1X_{P'}=d_1X_{Q'}$, (where $d_1$ corresponds to 
the operation of constructing edge decks of $n$-vertex 
and $(N-m+1)$-edge graphs): 
this follows from the fact that  
for any $F$, $e\in E(F)$ and $f\not \in E(F)-{e}$, 
$(F-e+f)^c = (F-e)^c -f$.
Now, if $m > N/2$  then $X_P=X_Q$, therefore, $d_1X_P=d_1X_Q$.
Else, $N-m+1 > N/2$, implying $X_{P'}=X_{Q'}$, 
therefore, $d_1X_P=d_1X_Q$.

\begin{cor}
Given the modified $1$-deck of a graph, its $1$-edge deck 
can be constructed.
\end{cor}

\subsection*{Remarks } 
\begin{enumerate} 
\item In his 1977 survey, Bondy asked if Lov\'asz's result could 
be proved constructively, (see problem 12 in  [B]). 
The inversion of the equation in Lemma~\ref{DelD} gives 
\[ D_k= \sum_{i=0}^k (-1)^{k+i}\ncr{m-i}{k-i}\Delta_i
\]
Therefore, when $m > N/2$, $D_m$ is zero.
Therefore,
\[
\Delta_0 = \Delta_1 - \Delta_2 + \Delta_3 - ... +(-1)^{m+1}\Delta_m
\]
and,
\[
X_P = \Delta_0X_P = \{\Delta_1 - \Delta_2 + \Delta_3 - ... 
+(-1)^{m+1}\Delta_m\}X_P  
\]
for any collection $P$ of graphs.  

Thus we explicitly know the reconstructed graph or a collection of graphs.

\item We can also prove Lov\'asz's result directly from $\Delta_1$ using
some properties of the Johnson graph defined below.

\begin{defn}
Johnson graph is a simple graph whose vertex set is 
the family of $m$-sets
of an $N$-set. Two vertices $U$ and $V$ are adjacent if and only 
if $|U \cap V|=m-1$.
\end{defn}
 
Let $J$ be the adjacency matrix of the Johnson graph with parameters 
$N=\ncr{n}{2}$ and $m$.  
Let the square matrix $B$ be defined as follows. The rows and columns 
of $B$ are indexed by all the labelled $m$-edge graphs on a fixed set 
of $n$ vertices, and $ij$-th entry is the number of ways of removing an 
edge from $G_j$ and adding an edge to get $G_i$. Note that the 
diagonal entry is $m$, since we can add the same edge that is removed.
Other entries of $B$ are either $0$ or $1$.
The matrix $A$ is defined similarly for unlabelled graphs with $m$ edges and
$n$ vertices. Thus matrix $A$ is the matrix $\Delta_1$ 
of the previous section. Matrix $P$ is defined by indexing the rows by 
unlabelled graphs and columns by labelled graphs, and the $ij$-th 
entry is $1$ if 
the labelled graph $G_j$ is isomorphic to the unlabelled graph $G_i$.  
Other entries are $0$.
As in [ER], one can verify that $AP=PB$, and every eigenvalue of 
$A$ is also an eigenvalue of $B$. But $B=mI+J$, and $-m$ is 
an eigenvalue of the Johnson graph if and only if $m\leq N/2$, 
(see [BCN]).  
So all eigenvalues of $A$ are nonzero when $m >  N/2$. This implies 
Lov\'asz's result for any collection of graphs. 
Theorem 2.7 can then be 
proved from this, but we have included the first proof because it is 
constructive. 

Perhaps even the $k$-edge version can be proved this way. 
In case of 
$k$-edge problem,
$B$ is simply $\sum_{i=0}^k\ncr{m-i}{k-i}J_i$ where $J_0=I$ and $J_i$
is the matrix of the $i$-th relation of Johnson scheme. 
The eigenvalues of matrices in a Johnson algebra are known: the $j$-th 
eigenvalue ($j\leq \mbox{min}(m,N-m)$) of $B$ is 
\[
\sum_{i=0}^k \sum_{l=0}^i(-1)^l\ncr{j}{l}\ncr{m-j}{i-l}\ncr{N-m-j}{i-l}\ncr{m-i}{k-i}
\] 
We need to prove that when $2m-k+1> N$, above polynomial in $j$
doesn't vanish in the range of integers $0$ to $\mbox{min}(m,N-m)$.

\item We have not solved the $k$-edge version of the problem of 
Krasikov and Roditty. But, 
since the procedure in the proof of Theorem~\ref{Deld} 
can be used also in the $k$-edge case, it 
it would be enough to prove 
the $k$-edge deck version 
of Lov\'asz's result.  

To prove the $k$-edge deck version of Lov\'asz's result 
we have indicated one approach using Johnson schemes. Another 
approach would be to relate $\Delta_i$, $i\geq k$ 
to $\Delta_k$ explicitly, and then prove Lov\'asz's 
result using those decks.  

\end{enumerate} 

\section*{Acknowledgements} 
This work was done at Combinatorics and Optimization, University of Waterloo,
Canada, and supported by the grant A7331 from NSERC. 
I take this opportunity to thank Adrian Bondy who arranged my visit.
He and Hugh Hind were very helpful and encouraging througout my stay at
Waterloo. I would also like to thank the referees for their many
valuable suggestions.
 
\section*{References}
\begin{itemize}
\item[{[B]}] J.A. Bondy  and R. L. Hemminger, Graph reconstruction - a survey,
{\em J. Graph Theory} {\bf 1} no. 3 (1977) 227-268. 
\item[{[BCN]}] A. E. Brouwer, A. M. Cohen and A. Neumaier,
{\em Distance-Regular Graphs}, Springer, Berlin, 1989.
\item[{[ER]}]  M. N. Ellingham and G. F. Royle, Vertex-switching
reconstruction of subgraph numbers and triangle-free graphs, 
{\em J. Combinatorial 
Theory Ser. B}{\bf  54} no. 2 (1992) 167-177.
\item[{[KR]}] I. Krasikov and Y. Roditty, Balance equations for reconstruction
problems, {\em Arch. Math. (Basel)} {\bf 48} no. 5 (1987) 458-464.
\item[{[L]}]  L. Lov\'asz, A note on the line reconstruction problem, {\em 
J. Combinatorial 
Theory Ser. B} {\bf  13} (1972) 309-310.
\item[{[M]}]  V. M\"uller, The edge reconstruction hypothesis is true for
graphs with more than $n.\log_2 n$ edges, {\em  J. Combinatorial Theory Ser. B}
{\bf 22} no. 3 (1977) 281-283.
\end{itemize}
\end{document}